# Complex Laplacian-based Distributed Control for Multi-Agent Network *

Aniket Deshpande     Pushpak Jagtap     Prashant bansode
Arun D. Mahindrakar     Navdeep M. Singh

**Abstract**

The work done in this paper, proposes a complex Laplacian-based distributed control scheme for convergence in the multi-agent network. The proposed scheme has been designated as cascade formulation. The proposed technique exploits the traditional method of organizing large scattered networks into smaller interconnected clusters to optimize information flow within the network. The complex Laplacian-based approach results in a hierarchical structure, with formation of a meta-cluster leading other clusters in the network. The proposed formulation enables flexibility to constrain the eigen spectra of the overall closed-loop dynamics, ensuring desired convergence rate and control input intensity. The sufficient conditions ensuring globally stable formation for proposed formulation are also asserted. Robustness of the proposed formulation to uncertainties like loss in communication links and actuator failure has also been discussed. The effectiveness of the proposed approach is illustrated by simulating a finitely large network of thirty vehicles.

## 1 Introduction

The broad application domain of distributed control methods in multi-agent systems have attracted a considerable attention of researchers in recent years. Some major application areas include formation control of unmanned air vehicles (UAV) [1], cooperative control of mobile vehicles [2], and distributed sensor networks [3]. The key challenges include minimization of control efforts, improvization in convergence time, communication constraints, and reducing computational cost.

Formation control problems have been serving as benchmark problems for decentralized control of multi-agent systems [4, 5]. Significant amount of research has been done to develop methods for addressing issues related to consensus problems such as cooperative control of multi-agnt systems subjected

---

to switching networks and communication delay [6, 7], collision avoidance [8], robustness to link and node failure [9], time varying formation control [10], and consensus of agent under input saturation [11]. Different techniques and frameworks have also been proposed to address agreement problems such as graph Laplacian-based method for formation stabilization [12], planar formation using complex Laplacian [13], formation control of heterogeneous nonlinear agent using passivity framework [14], leader follower architecture [15, 16], and hierarchical formation control [17].

The recent work in [13], has reported a complex Laplacian-based approach to achieve rigid planar formation. The work provides algebraic and geometrical conditions that ensure a globally stable formation. Laplacian-based methods often result in a slower moving average, as the number of communicating neighbours increase, delaying the network consensus which in turn, affects the system adversely. Reported literature intends to present a formulation to address such issues.

For this purpose, the conventional clustering method is exploited to systematically reorganize a large complex network into a number of distributed clusters. Such an arrangement is expected to channelize information exchange within the network leading to faster convergence. The inclusion of complex Laplacian-based methodology results in a hierarchical structure of the clusters led by a meta-cluster. The orientation and formation of clusters are controlled by the meta-cluster that commands the co-leaders of clusters. The method incorporates a bidirectional 2- rooted graph topology [18] that provides additional flexibility to control the overall formation along four degrees of freedom (viz, translation, rotation, and scaling) through the co-leaders of the meta-cluster.

Rest of the paper is organized as follows. Section II describes the preliminaries of graph theory and complex Laplacian with its necessary and sufficient conditions. Section III, includes discussion on the *cascade formulation* method to shape the information flow in large multi-agent systems by divides a large network into decoupled stable clusters with local decentralized control law. Section IV summarizes simulation results and comparative analysis of proposed formulation. The conclusions and open problems are discussed in Section V.

## 2 Preliminaries

The symbol $\mathbb{N}_k$ denotes the set of natural numbers greater than $k-1$. We use $\mathbb{C}^{n \times m}$ to denote a vector space of complex valued matrices with $n$ rows and $m$ columns. For $c \in \mathbb{C}$, $\text{Re}(c)$ and $\text{Im}(c)$ represent the real and imaginary parts of a complex number $c$, respectively. For a square matrix $F \in \mathbb{C}^{n \times n}$, $eig(F) = \lambda_1, \lambda_2, \ldots, \lambda_n \in \mathbb{C}$ represents eigenvalues of $F$ and the largest eigenvalue of $F$ is given as $\lambda_{\max}(F) = \max\{\text{Re}(\lambda_1), \text{Re}(\lambda_2), \ldots, \text{Re}(\lambda_n)\}$.

The interaction topology in a multi-agent system is represented using bidirectional graph $\mathcal{G} = (\mathcal{V}, \mathcal{E})$ with $n$ nodes $\mathcal{V} = \{1, 2, \ldots, n\}$ and edges $\mathcal{E} \subseteq \mathcal{V} \times \mathcal{V}$. Let the neighbor set of $i^{th}$ agent be defined as $\mathcal{N}_i = \{j \in \mathcal{V} | (j, i) \in \mathcal{E}\}$, where $i \in \{1, 2, \ldots, n\}$.

In order to prove results of this paper and to select proper interaction topology, it is important to introduce two definitions from [18].

**Definition 2.1** *For a bidirectional graph $\mathcal{G}$, a node $v \in \mathcal{V}$ is said to be 2-reachable from a non-singleton set $\mathcal{U}$ of nodes, if it is possible to reach node $v$ from any node in $\mathcal{U}$ after eliminating any one node except node $v$.*

**Definition 2.2** *A bidirectional graph $\mathcal{G}$ is said to be 2-rooted, if there exists a subset of two nodes, from which every other node is 2-reachable. These two nodes are termed as roots of the graph $\mathcal{G}$.*

Readers can refer to [13], for the graphical explanation of these definitions. The complex Laplacian $L$ for bidirectional graph $\mathcal{G}$ is given as

$$L(i,j) = \begin{cases} -w_{ij}, & \text{if } i \neq j \text{ and } j \in \mathcal{N}_i, \\ 0, & \text{if } i \neq j \text{ and } j \notin \mathcal{N}_i, \\ \sum_{j \in \mathcal{N}_i} w_{ij}, & \text{if } i = j, \end{cases} \quad (1)$$

where $w_{ij} \in \mathbb{C}$ is the complex weight associated with edge $(i,j)$. The definition of complex Laplacian also ensures that the row sum should be equal to zero (*i.e.,* it has at least one eigenvalue at origin with the corresponding eigenvector $\mathbf{1}_n$).

## 3 Planar formation using complex Laplacian

Consider a group of $n$ agents in a plane, with an objective to achieve a desired formation using distributed control laws. The control laws are assumed to be implementable with local information like relative distances with neighbors. In complex Laplacian approach, the formation configuration or shape of final formation is represented by assigning location $\xi_i \in \mathbb{C}$ in the complex plane to $i^{th}$ agent of the group. This complex formation vector $\xi = [\xi_1, \xi_2, \ldots, \xi_n]^T \in \mathbb{C}^n$ is referred to as *formation basis*. The $F_\xi$ is a function that acts on the formation basis $\xi$ along four degree-of-freedoms (translation, rotation, and scaling) to steer the formation as per the requirement and is written as,

$$F_\xi = c_1 \mathbf{1}_n + c_2 \xi, \quad c_1, c_2 \in \mathbb{C}.$$

Here, the bidirectional graph $\mathcal{G}$ with $n$ nodes is used as a *sensing graph* with edges $(i,j)$ representing a measure of relative position between agent $j$ and agent $i$ as $(z_j - z_i)$, where $z_j, z_i \in \mathbb{C}$ denotes the positions of $j^{th}$ and $i^{th}$ agents, respectively.

Suppose each agent is modeled as a fully actuated point mass with single-integrator kinematics given by,

$$\dot{z}_i = u_i, \quad i \in \{1, 2, \ldots, n\} \quad (2)$$

where $u_i \in \mathbb{C}$ is the velocity control input and the saturation limits are considered as $v_{\min} \leq \text{Re}(u_i) \leq v_{\max}$, $v_{\min} \leq \text{Im}(u_i) \leq v_{\max}$, $v_{\min}$ and $v_{\max}$ are the minimum and maximum velocities respectively, and are dependent on the actuator saturation limit. The local distributed control law to achieve a stable formation is then written as,

$$u_i = d_i \sum_{j \in \mathcal{N}_i} w_{ij}(z_j - z_i), \quad i \in \{1, 2, \ldots, n\} \tag{3}$$

where $d_i \in \mathbb{C}$ is a design parameter which decides the performance and global stability of the formation and $w_{ij}$ is the complex weight on edge $(i, j)$ represented in complex Laplacian $L$. The overall dynamics of the $n$ agent system with control law (3) is

$$\dot{z} = -DLz, \tag{4}$$

where $z = [z_1, z_2, \ldots, z_n]^T \in \mathbb{C}^n$, $L$ is the complex Laplacian of the communication graph $\mathcal{G}$, and $D = diag(d_1, d_2, \ldots, d_n)$ is stabilizing diagonal matrix. The diagonal matrix $D$ transforms the eigenvalues of (4) to the left half of the complex plane. The necessary and sufficient conditions to a design complex Laplacian $L$ and stabilizing matrix $D$ are discussed next.

## 3.1 Necessary and sufficient conditions

The necessary and sufficient conditions for construction of complex Laplacian and stabilizing matrix are stated in the following Lemmas [13].

**Lemma 3.1** *Let the formation basis $\xi \in \mathbb{C}^n$ satisfy $\xi_i \neq \xi_j$, $\forall i, j$. The equilibrium state of (4) forms a globally stable geometric formation $F_\xi$ if and only if there exists matrices $D, L \in \mathbb{C}^{n \times n}$ satisfying $eig(-DL) \leq 0$, $L\xi = 0$, and $rank(L) = n - 2$.*

**Lemma 3.2** *For a bidirectional graph $\mathcal{G}$ and $\xi \in \mathbb{C}^n$ satisfying $\xi_i \neq \xi_j$, $\forall i, j$. The algebraic conditions, $rank(L) = n-2$. and $L\xi = 0$ satisfies for all $L \in \mathbb{C}^{n \times n}$ if and only if $\mathcal{G}$ is 2-rooted.*

The Lemma 3.1 presents necessary algebraic condition to guarantee stationary formation and Lemma 3.2 gives graphical sufficiency condition to satisfy necessary condition mentioned in Lemma 3.1. The stability of the closed-loop system (4) depends on eigenvalues of $DL$. As complex Laplacian $L$ has its eigenvalues scattered over the whole complex plane, it is not always stable. Thus, it is important to design a diagonal matrix $D$ which stabilizes the closed-loop dynamics.

**Lemma 3.3** *If the bidirectional graph $\mathcal{G}$ is 2-rooted then there exists a stabilizing matrix $D$ for the system $\dot{z} = -Lz$ such that (4) is stable.*

Proof's for the above lemmas are not included for the obvious and readers are referred to [13] for the same.

## 3.2 Performance analysis of consensus algorithm

Eigenvalues of the stabilized complex Laplacian matrix exhibit important information about performance of consensus algorithm such as stability, convergence rate, and control efforts. Some properties of eigenvalues of complex Laplacian are

*1.* Complex Laplacian $L$ for 2-rooted graph topology has two eigenvalues at origin with corresponding eigenvectors $\mathbf{1}_n$ and formation basis $\xi$.

*2.* Unlike a real-valued Laplacian, complex Laplacian may have eigenvalues in left half of complex plane.

*3.* All non-zero eigenvalues of $L$ can be shifted to right half of complex plane by pre-multiplying it with real valued invertible diagonal matrix $D$ without affecting property 1 [19].

Let the stabilized complex Laplacian $DL$ have its eigenvalues at $\lambda_1 = \lambda_2 = 0 < \text{Re}(\lambda_3) \leq \text{Re}(\lambda_4) \leq \cdots \leq \text{Re}(\lambda_n)$. The smallest non-zero eigenvalue of stabilized complex Laplacian matrix $\lambda_3$ is considered as an extension to the concept of algebraic connectivity [20] of real-valued Laplacian to complex Laplacian with 2-rooted graph topology and is used as a measure of performance of collective dynamics. The largest eigenvalue of $DL$, $\lambda_n$, is used as a measure of the intensity of control signal [21]. It is very important to limit control signal magnitude (*i.e.,* by placing $\lambda_n$ properly), as it can produce instability due to saturation. Our main objective is to strengthen the algebraic connectivity of the complex Laplacian while limiting intensity of the control inputs to improve the speed of convergence in large multi-agent system. This is achieved by converting consensus problem into special formulation by cascading small clusters as discussed in the next section.

# 4 Proposed Formulation

The complex Laplacian-based consensus algorithm requires 2-rooted graph topology with two roots acting as co-leaders for orientation and scaling of formation. It is observed that as the number of agents increases, $\lambda_{\max}(DL)$ increases which may result in instability due to saturation of control inputs [22, 23]. The system can be stabilized by scaling down the complex Laplacian by appropriate factor $k$ as $kDL$, where $k \in (0,1)$. However, it affects the algebraic connectivity of network which in turn affects the convergence time. This sets the trade-off between convergence time and the control input intensity rendering the design of stabilizing matrix $D$, a tedious task. Another major issue with complex Laplacian-based control law is its inability to deal with communication and actuation failure of an agent. The overall dynamics lead to instability in case of failure of communication link or actuator of an agent due to interruption in information flow. Robustness can be incorporated in a complex Laplacian-based control by implementing the proposed systematic methodology of *Cascade Formulation*.

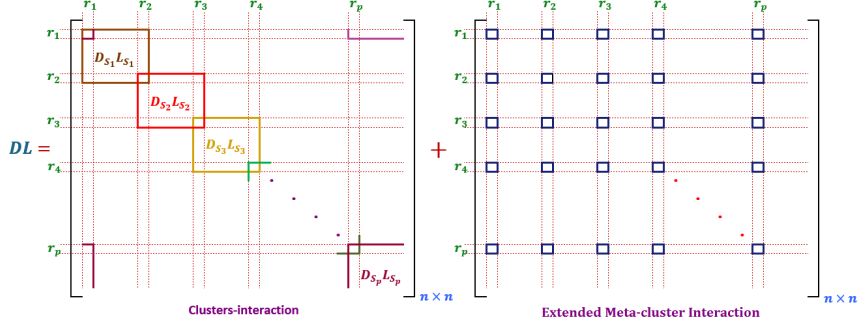

Figure 1: Representation of complex Laplacian structure using cascade formulation

## 4.1 Methodology

Consider a mutli-agent system with $n$-agents represented by bidirectional graph $\mathcal{G} = (\mathcal{V}, \mathcal{E})$ with $n$ nodes $\mathcal{V} = \{1, 2, \ldots, n\}$, $n \in \mathbb{N}_3$ and edges $\mathcal{E} \subseteq \mathcal{V} \times \mathcal{V}$. A systematic approach designated as, *cascade formulation* will be elaborated to establish inter-agent interconnection and overcome aforementioned issues. The proposed formulation divides the large multi-agent system with $n$ agents, into a hierarchical structure of $p$ clusters led by a meta-cluster (defined in Definition 4.1). The $p$ clusters are denoted by $S_i^{q_i \times q_i}$, where $i \in \{1, 2, \cdots, p\}$ and $q_i \in \mathbb{N}_3$ is a number of agents in each cluster. Every cluster $S_i^{q_i \times q_i}$ satisfies the properties of 2-rooted graph topology (*i.e.*, every cluster has $q_i - 2$ follower agents and two roots which act as co-leaders of the respective cluster). The root $r_i$ of every cluster is shared between its adjacent clusters as shown in Figure 1. Every cluster satisfies algebraic and geometric conditions given in Lemma 3.1 and Lemma 3.2 which states that every cluster has stabilized complex Laplacian $D_{S_i} L_{S_i}$ associated with it. Now let us introduce a meta-cluster as follows,

**Definition 4.1** *A set of nodes, $M \in \mathcal{V}$, in cascade formulation is said to be meta-cluster, if one has*

$$M = \Big\{r_i \big|\ r_i, r_j\ \in\ \mathcal{V} \text{ and } (i,j) \in \mathcal{E},\ \forall\ i \neq j;\ i, j \in \{1, 2, \cdots, p\}\Big\},$$

*where $r_i$ is the root of a cluster provided, all nodes in $M$ are connected by a bidirectional $2$-rooted graph.*

The roots of meta-cluster act as main co-leaders of the network (overall formation). The orientation and scaling of the network is then decided by the formation basis of meta-cluster $\xi_M$. Moreover, $\xi_M$ is the vector of cluster agents in complex coordinates. The stabilized complex Laplacian for meta-cluster is denoted by $D_M L_M$.

**Lemma 4.2** *Consider a multi-agent network of n-agents having closed loop dynamics $\dot{z} = -DLz$, interconnected in cascade formulation. Then the eigenvalues of each cluster $S_i$ and meta-cluster $M$ are independent of each other if each cluster and the meta-cluster satisfy 2-rooted graph topology and the clusters are connected only through roots $r_i$.*

**Proof: 4.3** *Consider the stabilized complex Laplacian is designed for each of the clusters and meta-cluster independently as $D_{S_i}L_{S_i}$ and $D_M L_M$, respectively. This means that the clusters and meta-cluster satisfy conditions in Lemma 3.1 and Lemma 3.2. By using similarity transform of matrix to transform stabilized complex Laplacian into diagonal matrix with diagonal entries as corresponding eigenvalues [24], [25].*

*Let $P_{S_i} \in \mathbb{C}^{q_i \times q_i}$ and $P_M \in \mathbb{C}^{p \times p}$ be the matrices of right eigenvectors of the cluster $S_i$, $i \in \{1, 2, \ldots, p\}$ and meta-cluster $M$, respectively. The corresponding diagonal matrices can be represented as*

$$\Lambda_{S_i} = P_{S_i}^{-1}(D_{S_i}L_{S_i})P_{S_i} = diag(0, \lambda_2^{S_i}, \lambda_3^{S_i}, \ldots, \lambda_{q_i-1}^{S_i}, 0), \tag{5}$$

$$\Lambda_M = P_M^{-1}(D_M L_M)P_M = diag(0, \lambda_2^M, \lambda_3^M, \ldots, \lambda_{p-1}^M, 0). \tag{6}$$

*The 2-rooted graph topology of $p$ clusters led by a meta-cluster ensures that there exist two eigenvalues at origin corresponding to roots of the graph. The first and last row of $\Lambda_{S_i}$ and $\Lambda_M$ in (5) and (6) represents roots of clusters and meta-cluster, respectively. The proposed formulation considers that the adjacent clusters are connected only through its roots as shown in Fig.1. This implies that eigenvalues of cluster $S_i$ will not affect the eigenvalues of its adjacent clusters.*

*The 2-rooted structure of meta-cluster and Definition 4.1 ensure that meta-cluster has $p-2$ non-zero eigenvalues corresponding to roots of clusters $r_i$. This results in the two zero roots acting as co-leaders of the network without affecting the other eigenvalues of clusters.*

**Remark 4.4** *Lemma 4.2 implies that each cluster and meta-cluster can be treated as decoupled systems. Therefore, it is possible to design stabilized complex Laplacian DL for individual cluster and meta-cluster.*

**Theorem 4.5** *Consider a multi-agent network of n-agents and Lemma 4.2. The closed-loop dynamics $\dot{z} = -DLz$ under cascade formulation results in a globally stable formation if individual cluster and meta-cluster are stabilized using complex Laplacian-based control law.*

**Proof: 4.6** *The proof of this theorem is a direct consequence of Lemma 4.2. Assume that the individual stabilizing diagonal matrix $D_{S_i}$ and $D_M$ for every cluster and meta-cluster are designed to place eigenvalues of $D_{S_i}L_{S_i}$ and $D_M L_M$ to right-hand side, respectively. As discussed in Subsection 4.1, each cluster and meta-cluster satisfies algebraic and geometric conditions.*

*Consider a cluster $S_i$ with its corresponding formation basis $\xi_{S_i}$, the algebraic condition is*

$$(D_{S_i}L_{S_i})\xi_{S_i} = 0, \quad \text{for } i \in \{1, 2, \ldots, p\} \tag{7}$$

and the algebraic condition for meta-cluster $M$ and formation basis $\xi_M$ is

$$(D_M L_M)\xi_M = 0. \tag{8}$$

At roots of a cluster, there is an interaction between its adjacent clusters and meta-cluster which gives,

$$\sum_{i \in \Omega_{r_j}} (D_{S_i} L_{S_i}^{r_j})\xi_{S_i} + (D_M L_M^{r_j})\xi_M = 0, \text{ for } j \in \{1, 2, \ldots, p\}, \tag{9}$$

where $\Omega_{r_j}$ is the set of adjacent clusters of root $r_j$, $D_{S_i} L_{S_i}^{r_j}$ and $D_M L_M^{r_j}$ represent row corresponding to $r_j^{th}$ root of cluster $S_i$ and meta-cluster $M$ respectively. By using (7), (8) and (9), one has algebraic condition of overall formulation for the complete formation basis $\xi$ as

$$(DL)\xi = 0. \tag{10}$$

In the cluster interaction (see Figure 1), there are $p$ zero eigenvalues corresponding to roots of clusters $r_1, r_2, \ldots,$ and $r_p$. This results in rank of cluster interaction as $n - p$ and the meta-cluster has $n - 2$ non-zero eigenvalues at roots because of 2-rooted topology. Thus the overall rank of proposed formulation is $n - 2$. This satisfies the algebraic condition given in Lemma 3.1 to achieve globally stable formation.

Let the algebraic connectivity and largest eigenvalue of $p$ clusters be denoted by $\lambda_a(D_{S_1} L_{S_1}), \lambda_a(D_{S_2} L_{S_2}), \ldots, \lambda_a(D_{S_p} L_{S_p})$ and $\lambda_{\max}(D_{S_1} L_{S_1}), \lambda_{\max}(D_{S_2} L_{S_2})$, $\ldots, \lambda_{\max}(D_{S_p} L_{S_p})$, respectively. For meta-cluster, it is represented as $\lambda_a(D_M L_M)$ and $\lambda_{\max}(D_M L_M)$, respectively. As mention in Lemma 4.2, the eigenvalues of each cluster is independent of other. Thus, the algebraic connectivity and largest eigenvalue of formulation are given as

$$\lambda_a = \min\{\lambda_a(D_{S_1} L_{S_1}), \lambda_a(D_{S_2} L_{S_2}), \ldots, \lambda_a(D_{S_p} L_{S_p}), \lambda_a(D_M L_M)\},$$
$$\lambda_{\max} = \max\{\lambda_{\max}(D_{S_1} L_{S_1}), \lambda_{\max}(D_{S_2} L_{S_2}), \ldots, \lambda_{\max}(D_{S_p} L_{S_p}), \lambda_{\max}(D_M L_M)\}.$$

**Remark 4.7** *One can easily increase the hierarchy of formulation by cascading several meta-clusters whose formation basis is described by meta-meta-cluster without affecting the stability and performance of overall formation.*

**Proposition 4.8** *The stabilized formation of multi-agent network of $n$-agents under proposed formulation, does not loose its overall stability even if it is subjected to uncertainties like actuation or communication link failure.*

**Proof: 4.9** *The proposition is a direct consequence of Lemma 4.2, Remark 4.4, and Theorem 4.5. We have seen that the stability achieved in proposed formulation is cluster-wise independent. Thus, even if one of the agent from a particular cluster losses its communication or actuation, it cannot affect the stability (i.e., eigenvalues) of adjacent clusters or meta-cluster. Thus a major part of the large network remains undisturbed and stable.*

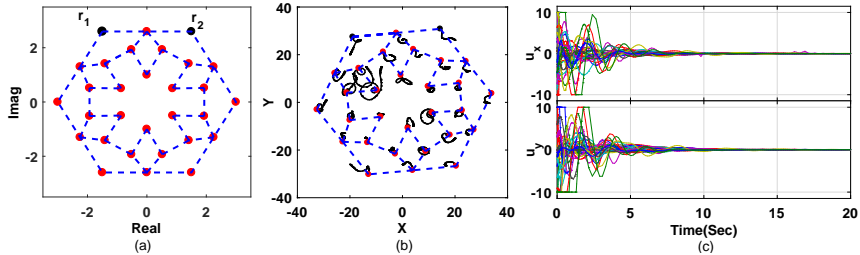

Figure 2: Using conventional approach: (a) Interconnection, (b) Closed-loop response, (c) Control inputs $u_x$ and $u_y$.

### 4.2 Designing stabilization matrices

Unlike real-valued Laplacian, the eigenvalues of a complex Laplacian are scattered randomly in the entire complex plane due to complex parameter matrix (which, in most cases do not appear in conjugate pairs). To the best of author's knowledge, the existing conventional optimization frameworks are not suitable to handle such complex parameter systems. Thus, an evolutionary algorithm based technique has been adopted to design stabilizing matrices. In the proposed formulation, the matrices $D_{S_i}$ and $D_M$ are designed using genetic algorithm [26, 27] with an objective of restricting the eigenvalue spectrum bandwidth of the complex Laplacian to a desired range. Note that, the genetic algorithm is only used to find one of the infinite solutions of stabilizing matrices ensuring desired range over eigenvalue spectrum.

The trade-off between convergence time and control efforts can be formulated in the form of an objective function as described below,

$$\min \begin{cases} \tau = |2\min(Re(eig(DL))) - \overline{\lambda}_{\min} - \overline{\lambda}_{\max}|, \\ \sigma = |2\max(Re(eig(DL))) - \overline{\lambda}_{\min} - \overline{\lambda}_{\max}|, \end{cases} \quad (11)$$

where $\tau$ represents an objective function for the rate of convergence, $\sigma$ denotes the objective function for control input intensity, $\overline{\lambda}_{\min}$ and $\overline{\lambda}_{\max}$ are lower and upper bounds on the spectrum of nonzero eigenvalues of complex Laplacian.

The objective function $\tau$ ensures that the all non-zero eigenvalues of clusters and meta-cluster are greater than required algebraic connectivity, $\lambda_a$, which controls the rate of convergence. The magnitude of control input $u_i$ is restricted by bounding the eigenvalues below $\overline{\lambda}_{\max}$ using objective function $\sigma$. This helps to avoid instability in case of saturation on control inputs. One can select the value of $\overline{\lambda}_{\max} > \lambda_a$ and close to the desired algebraic connectivity $\lambda_a$.

## 5 Simulation and Results

In this section, we present comparative simulation results to show efficacy and robustness of proposed formulation.

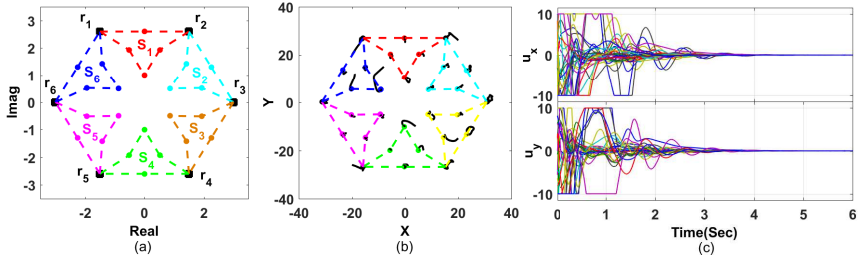

Figure 3: Using cascade formulation: (a) Clusters, meta-cluster and their interaction, (b) Closed-loop response, (c) Control inputs $u_x$ and $u_y$.

## 5.1 Performance analysis

For the purpose of simulation, a network of 30 agents modeled by single integrator kinematics as given in (2) has been considered. The velocity input constraints are $-10 \leq \text{Re}(u_i) \leq +10$ and $-10 \leq \text{Im}(u_i) \leq +10$. The system is simulated using MATLAB® Simulink. The target formation of agents represented by formation basis $\xi$ in complex plane and communication topology is shown in Figure 2(a). The black nodes indicate roots of the graph $\mathcal{G}$. Simulation trajectories using control law (3) are shown in Figure 2(b). The stabilized complex Laplacian has its algebraic connectivity at $\lambda_3 = 0.0027 + \imath 0.1447$ and the largest eigenvalue at $\lambda_{\max} = 16.44 - \imath 6.6786$. Convergence time and control signals are shown in Figure 2(c).

The proposed *cascaded formulation* is applied to the network discussed above by dividing the network into six homogeneous clusters $(S_1, S_2, \ldots, S_6)$ as shown in Figure 3(a). The black nodes represent roots of clusters that comprise a meta-cluster. Considering uniformity in clusters, a single complex Laplacian $L_S$ and a stabilizing matrix $D_S$ can be designed for all clusters. The algebraic connectivity and largest eigenvalue of overall formation are at $1.5762 - \imath 3.4779$ and $23.352 - \imath 2.4619$, respectively. The response using proposed approach is given in Figure 3(b). Figure 3(c) shows control signals $u_x$ and $u_y$. It is observed that the structured and distributed information flow due to the proposed algorithm reduces the convergence time while satisfying constraints on the control inputs.

## 5.2 Robustness to communication and actuation failure

We now illustrate the robustness property of proposed formulation for the case of failure in any communication link or actuator of an agent. It is easily observed from stability conditions that the conventional complex Laplacian-based control law causes instability due to communication interruption and actuation failure.

The proposed formulation incorporates stability in individual clusters and the meta-cluster. Due to this, the communication and actuation failure of any agent of a particular cluster will not affect the stability of other clusters. Results have been illustrated considering two cases. In the first case, failure of communication link of an agent within cluster is considered by making its corresponding edge weight in the Laplacian equal to zero. Simulation result is shown in Figure 4(a). The result establishes the fact that the failure will not affect the stability of other clusters and meta-cluster. In the second case, failure of link in meta-cluster agent is simulated similarly and it is observed that all clusters approach stable formation. The failure only affects the orientation and scaling of clusters adjacent to affected link (see Figure 4(b)).

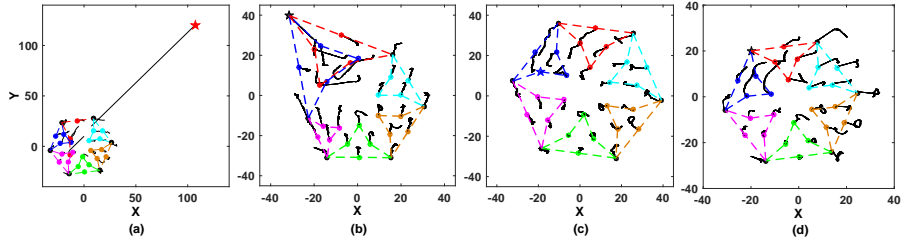

Figure 4: Response to loss of communication in (a) cluster agent, (b) meta-cluster agent; Response to actuator failure in (c) cluster agent, (d) meta-cluster agent, where star indicates agent with failure

The simulation results for actuation failure of an agent in cluster and meta-cluster are shown in Figures 4(c) and 4(d), respectively. One can readily observe that the proposed formulation forms stable formation around the agent even under actuation failure.

# 6  Conclusions and Open Problems

A novel approach is formulated to solve formation control problem in large multi-agent systems while attaining robustness to communication link and actuation failure. The cascade formulation proposed in this paper channelizes information flow throughout the network efficiently. The formulation divides the complex network into small clusters to incorporate decentralized information exchange between meta-cluster and agents of individual clusters. The 2–rooted bidirectional graph topology is adopted to form clusters and meta-cluster which allows them to be a decoupled dynamical systems. This offers flexibility in designing individual control laws, that satisfies the bounds on control inputs and achieves stable formation in desired convergence time. Moreover, it is illustrated that the proposed formulation is relatively robust even if the information flow in network is subjected to uncertainties like communication link and actuation failure in agents. The cascade formulation also provides for organization of distributed clusters at different hierarchies in complex systems, which is helpful in many applications like synchronization and collective task handling in multi-robot systems.

As an extension to the proposed methodology, our future research is focused on some challenging issues like collision avoidance, self organizing cascade formulation, etc. Moreover, the control strategy can be cast into well structured optimization framework to design a stabilizing diagonal matrix for optimization of control efforts and convergence time.

# References

[1] Reza Olfati-Saber. Flocking for multi-agent dynamic systems: Algorithms and theory. *IEEE Transactions on Automatic Control*, 51(3):401–420, 2006.

[2] Ali Jadbabaie, Jie Lin, et al. Coordination of groups of mobile autonomous agents using nearest neighbor rules. *IEEE Transactions on Automatic Control*, 48(6):988–1001, 2003.

[3] Petter Ögren, Edward Fiorelli, and Naomi Ehrich Leonard. Cooperative control of mobile sensor networks: Adaptive gradient climbing in a distributed environment. *IEEE Transactions on Automatic Control*, 49(8):1292–1302, 2004.

[4] Reza Olfati-Saber, Alex Fax, and Richard M Murray. Consensus and cooperation in networked multi-agent systems. *Proceedings of the IEEE*, 95(1):215–233, 2007.

[5] J Alexander Fax and Richard M Murray. Information flow and cooperative control of vehicle formations. *IEEE Transactions on Automatic Control*, 49(9):1465–1476, 2004.

[6] Reza Olfati-Saber and Richard M Murray. Consensus problems in networks of agents with switching topology and time-delays. *IEEE Transactions on Automatic Control*, 49(9):1520–1533, 2004.

[7] Shuai Liu, Lihua Xie, and Frank L Lewis. Synchronization of multi-agent systems with delayed control input information from neighbors. *Automatica*, 47(10):2152–2164, 2011.

[8] Dong Eui Chang, Shawn C Shadden, Jerrold E Marsden, and Reza Olfati-Saber. Collision avoidance for multiple agent systems. *Proceedings of IEEE Conference on Decision and Control*, 2003.

[9] Wente Zeng and Mo-Yuen Chow. Resilient distributed control in the presence of misbehaving agents in networked control systems. *IEEE Transactions on Cybernetics*, 44(11):2038–2049, 2014.

[10] Lara Brinón-Arranz, Alexandre Seuret, and Carlos Canudas-de Wit. Cooperative control design for time-varying formations of multi-agent systems. *IEEE Transactions on Automatic Control*, 59(8):2283–2288, 2014.

[11] Housheng Su, Michael ZQ Chen, James Lam, and Zongli Lin. Semi-global leader-following consensus of linear multi-agent systems with input saturation via low gain feedback. *IEEE Transactions on Circuits and Systems I: Regular Papers*, 60(7):1881–1889, 2013.

[12] Zhiyun Lin, Mireille Broucke, and Bruce Francis. Local control strategies for groups of mobile autonomous agents. *IEEE Transactions on Automatic Control*, 49(4):622–629, 2004.

[13] Zhiyun Lin, Lili Wang, Zhimin Han, and Minyue Fu. Distributed formation control of multi-agent systems using complex laplacian. *IEEE Transactions on Automatic Control*, 59(7):1765–1777, 2014.

[14] Murat Arcak. Passivity as a design tool for group coordination. *IEEE Transactions on Automatic Control*, 52(8):1380–1390, 2007.

[15] Florian Dörfler and Bruce Francis. Geometric analysis of the formation problem for autonomous robots. *IEEE Transactions on Automatic Control*, 55(10):2379–2384, 2010.

[16] Pushpak Jagtap, Aniket Deshpande, Navdeep M Singh, and Faruk Kazi. Complex laplacian based algorithm for output synchronization of multi-agent systems using internal model principle. *IEEE Conference on Control Applications (CCA)*, pages 1811–1816, 2015.

[17] Ji-Wook Kwon and Dongkyoung Chwa. Hierarchical formation control based on a vector field method for wheeled mobile robots. *IEEE Transactions on Robotics*, 28(6):1335–1345, 2012.

[18] Lili Wang, Zhimin Han, Zhiyun Lin, and Gangfeng Yan. Complex laplacian and pattern formation in multi-agent systems. *24th IEEE Chinese Control and Decision Conference (CCDC)*, pages 628–633, 2012.

[19] CS Ballantine. Stabilization by a diagonal matrix. *Proceedings of the American Mathematical Society*, 25(4):728–734, 1970.

[20] Miroslav Fiedler. Algebraic connectivity of graphs. *Czechoslovak mathematical journal*, 23(2):298–305, 1973.

[21] Daniel Y Abramovitch and Gene F Franklin. On the stability of adaptive pole-placement controllers with a saturating actuator. *IEEE Transactions on Automatic Control*, 35(3):303–306, 1990.

[22] Tingshu Hu and Zongli Lin. *Control systems with actuator saturation: analysis and design*. Springer Science & Business Media, 2001.

[23] Sophie Tarbouriech, Germain Garcia, and Adolf Hermann Glattfelder. *Advanced strategies in control systems with input and output constraints*. Springer, 2007.

[24] Gene H Golub and Charles F Van Loan. *Matrix computations*, volume 3. JHU Press, 2012.

[25] Roger A Horn and Charles R Johnson. *Matrix analysis*. Cambridge university press, 2012.

[26] Mandavilli Srinivas and Lalit M Patnaik. Genetic algorithms: A survey. *Computer*, 27(6):17–26, 1994.

[27] Kim-Fung Man, Kit-Sang Tang, and Sam Kwong. Genetic algorithms: concepts and applications. *IEEE Transactions on Industrial Electronics*, 43(5):519–534, 1996.